\magnification=1200
\input amstex
\documentstyle{amsppt}

\def\wt#1{\widetilde {#1}}

\def\SE {\Cal E}

\def\SH {\Cal H}

\def\SL {\Cal L}

\def\SO {\Cal O}

\def\SU {\Cal U}

\topmatter
\title Minimal rational curves in moduli spaces of stable bundles
 \endtitle
\author Xiaotao Sun   \endauthor
\address Institute of Mathematics, Chinese Academy of Sciences,
Beijing 100080, China \endaddress \email
xsun$\@$math08.math.ac.cn\endemail
\address Department of Mathematics, The University of
Hong Kong, Pokfulam Road, Hong Kong\endaddress \email
xsun$\@$maths.hku.hk\endemail
\thanks The work is supported by a grant of NFSC for outstanding
young researcher at contract number 10025103.
\endthanks
\endtopmatter
\document

\heading{Introduction}\endheading

Let $C$ be a smooth projective curve of genus $g\ge 2$ and $\SL$ be
a line bundle on $C$ of degree $d$. Assume that $r\ge 2$ is an integer 
coprime with $d$. Let $M:=\SU_C(r,\SL)$ be the 
moduli space of stable vector bundles on $C$ of 
rank $r$ and with the fixed determinant $\SL$. It is well-known that
$M$ is a smooth projective Fano variety with Picard
number $1$. For any projective curve in $M$, we can define
its degree with respect to the ample anti-canonical line bundle $-K_M$. 
A natural question raised by Jun-Muk Hwang
(see Question 1 in [Hw]) is to determine all rational curves of
$minimal$ $degree$ passing through a generic point of $M$.
In this short note, we prove the following theorem

\proclaim{Main Theorem} Assume that $g\ge 4$. Then any rational curve of $M$ 
passing 
through a generic
point has degree at least $2r$. It has degree $2r$ if and only if
it is a Hecke curve.\endproclaim

Our idea is simple and the proof is elementary. If $E$ is a vector bundle on $X=C\times\Bbb P^1$
that induces the morphism of $\Bbb P^1$ to $M$. Then a simple computation
shows that its degree equals to the second Chern class of $\SE nd(E)$.
If the restriction of $E$ to the generic fiber of ruled surface
$f: X\to C$ is semistable, then one sees easily that $c_2(\SE nd(E))$ is
at least $2r$, and it is $2r$ if and only if $c_2(E)=1$ (after tensoring
$E$ by suitable line bundle pulling back from $\Bbb P^1$). This will
force $E$ to be an extension
$$0@>>> f^*V@>i>> E@>\phi>>\SO_{\{p\}\times\Bbb P^1}(-1)@>>> 0,$$
where $V$ is a bundle on $C$. That is, after performing elementary
transformation on $E$ $once$ along $one$ fiber, $E$ becomes a pullback
of a vector bundle on $C$. For any $x\in\Bbb P^1$, restricting above sequence
to $C\times\{x\}$ and denote $E|_{C\times\{x\}}$ by $E_x$, we have
$$0@>>> V@>i_x>> E_x@>\phi_x>>\SO_{\{p\}\times\Bbb P^1}(-1)_x@>>> 0.$$
Let $\iota_x:V_p\to E_x|_p=E_{(p,x)}$ be the homomorphism between the fibers
at $p$ induced by the sheaf map $i_x$. Then the Hecke modifications
$\{({\wt W}^{ker(\iota_x)})^{\vee};\,\,x\in \Bbb P^1\}$ of $V$ along 
$\{ker(\iota_x)\subset V_p;\,\,x\in \Bbb P^1\}$ are exactly
$\{E_x;\,\,x\in\Bbb P^1\}$. Thus the given curve is a Hecke curve by
definition.

If the restriction of $E$ to the generic fiber is not semistable, then
using the so called relative Harder-Narasimhan filtration we are able
to prove that $c_2(\SE nd(E))>2r$. Here we also need the condition that
the rational curve passes through the generic point of $M$ which
corresponds a $(1,1)$-stable bundle when $g\ge 4$.    

\heading{\S1 Hecke curves }\endheading 

Since we are discussing the question raised in [Hw], we copy the definition
of Hecke curves in [Hw] and adopt notation there.
Given two nonnegative integers $k$, $\ell$, a vector bundle $W$ of rank $r$
and degree $d$ on $C$ is $(k,\ell)$-stable, if, for each proper subbundle
$W'$ of $W$, we have
$$\frac{deg(W')+k}{rk(W')}<\frac{deg(W)+k-\ell}{r}.$$
The usual stability is equivalent to $(0,0)$-stability. The dual bundle
of a $(k,\ell)$-stable bundle is $(\ell,k)$-stable.

\proclaim{Lemma 1.1 ([NR])} If $g\ge 4$, a generic point $[W]\in M$
corresponds to a $(1,1)$-stable bundle $W$.\endproclaim

\proclaim{Lemma 1.2 ([NR])} Let $0\to V\to W\to \SO_p\to 0$ be an exact sequence,
where $\SO_p$ is the $1$-dimensional skyscraper sheaf at $p\in C$.
If $W$ is $(k,\ell)$-stable, then $V$ is $(k,\ell-1)$-stable.
\endproclaim 

Let $[W]\in M$ be a generic point corresponding to a $(1,1)$-stable 
bundle $W$ over $C$. We will use $V^{\vee}$ to denote
the dual vector bundle (or dual vector space) of a vector bundle $V$ 
(or a vector space $V$). Let $\Bbb P(W)$ be the projective bundle consisting of lines through
the origin on each fiber. For $p\in C$ and $\zeta\in \Bbb P(W^{\vee}_p)$, 
define
a vector bundle $W^{\zeta}$ by
$$0\to W^{\zeta}\to W\to (W_p/\zeta^{\bot})\otimes\SO_p\to 0\tag 1.1$$
where $\zeta^{\bot}$ denotes the hyperplane in $W_p$ annihilated by 
$\zeta$. Let $\iota:W_p^{\zeta}\to W_p$ be the homomorphism between 
the fibers at $p$ induced by the sheaf injection $W^{\zeta}\to W$.
The kernel $ker(\iota)$ of $\iota$ is a $1$-dimensional subspace of
$W_p^{\zeta}$ and its annihilator $ker(\iota)^{\bot}$ is a hyperplane in
$(W^{\zeta})^{\vee}_p$. Let $\SH$ be a line in $\Bbb P(W^{\zeta}_p)$ containing the
point $[ker(\iota)]$. For each point $[l]\in\SH$ corresponding to a $1$-dimensional
subspace $l\subset W^{\zeta}_p$, define a vector bundle ${\wt W}^l$ by
$$0\to {\wt W}^l\to (W^{\zeta})^{\vee}\to 
((W^{\zeta})^{\vee}_p/l^{\bot}))\otimes\SO_p\to 0\tag1.2$$
where $l^{\bot}\subset(W^{\zeta})^{\vee}_p$ is the hyperplane annihilating 
$l$. This bundle ${\wt W}^l$ is stable for each $[l]\in \SH$ by
Lemma 1.2. It is easy to check that
for $l=ker(\iota)$, 
$${\wt W}^{ker(\iota)}\cong W^{\vee}.\tag 1.3$$

\newpage

Thus $\{({\wt W}^l)^{\vee}\,; \,l\in\SH\}$ defines a rational curve passing
through $[W]\in M$. A rational curve on $M$ constructed in this way is called a
Hecke curve. By using [NR], it can be shown that a Hecke curve is smooth
and has degree $2r$ with respect to $-K_M$.
There is an equivalent description for $({\wt W}^l)^{\vee}$.
For any vector bundle, say $W^{\zeta}$, on $C$ and a subspace
$K$ of the fiber $W^{\zeta}_p$ ($p\in C$).
According to [NS], there are two canonical constructions called
Hecke modifications (see Remark 2.4 of [NS]). For any $[l]\in\SH$,
the stable bundle $({\wt W}^l)^{\vee}$ is in fact obtained from
$W^{\zeta}$ by performing the second Hecke modification along
the subspace $l\subset W^{\zeta}_p$. 

\heading{\S2 Proof of the Main Theorem }\endheading 

For any rational curve
$\Bbb P^1\subset M$ through a general point of $M$,
let $E$ be the vector bundle on $X:=C\times\Bbb P^1$, which induces
the embedding $\Bbb P^1\subset M$.
Let $\pi:X=C\times\Bbb P^1\to \Bbb P^1$ be the projection and
$\Bbb E\subset\SE nd(E)$ be the subbundle of trace free. Then, since
$\pi_*(\Bbb E)=0$, we have $T_M|_{\Bbb P^1}=R^1\pi_*\Bbb E$ and,
 by using Leray spectral sequence and Riemann-Roch theorem,   
$$-\chi(\Bbb E)=\chi(R^1\pi_*\Bbb E)=-K_M\cdot\Bbb P^1+(r^2-1)(g-1).$$
By using $\chi(\Bbb E)=deg(ch(\Bbb E)\cdot td(T_X))_2$, noting 
$c_1(\Bbb E)=c_1(\SE nd(E))=0$, we get
$$-K_M\cdot\Bbb P^1=c_2(\Bbb E)=2rc_2(E)-(r-1)c_1(E)^2:=\Delta(E).\tag2.1$$

Let $f: X=C\times\Bbb P^1\to C$ be the projection. Then, for any torsion free
sheaf $E$ on the ruled surface $X$, its restriction to a generic fiber
$f^{-1}(\xi)=X_{\xi}$ has the form
$$E|_{X_{\xi}}=\bigoplus_{i=1}^n\SO_{X_{\xi}}(\alpha_i)^{\oplus r_i},\quad
\alpha_1>\,\,\cdots\,\,>\alpha_n.$$
The $\alpha=(\alpha_1^{\oplus r_1},\,...,\,\alpha_n^{\oplus r_n})$ is 
called the generic splitting type of $E$. In our case, tensoring $E$ by
$\pi^*\SO(-\alpha_n)$, we can (and we will) assume that $\alpha_n=0$.
Any such $E$ admits a relative Hardar-Narasimhan filtration
$$0=E_0\subset E_1\subset \,\cdots\,\subset E_n=E$$
of which the quotient sheaves $F_i=E_i/E_{i-1}$ are torsion free with 
generic splitting type $(\alpha_i^{\oplus r_i})$ respectively.
Then it is easy to see that
$$\aligned 2c_2(E)&=2\sum_{i=1}^nc_2(F_i)+2\sum^n_{i=1}c_1(E_{i-1})c_1(F_i)\\
&=2\sum_{i=1}^nc_2(F_i)+c_1(E)^2-\sum_{i=1}^nc_1(F_i)^2.\endaligned$$
Thus 
$\Delta(E)=
2r\sum_{i=1}^nc_2(F_i)+c_1(E)^2-r\sum_{i=1}^nc_1(F_i)^2.$

\newpage

Let $F_i'=F_i\otimes \pi^*\SO_{\Bbb P^1}(-\alpha_i)$ ($i=1,\,...,n$),
thus they have generic splitting type $(0^{\oplus r_i})$ respectively.
Let $c_1(F_i)=f^*\SO_C(d_i)+\pi^*\SO_{\Bbb P^1}(r_i\alpha_i)$, where
$\SO_C(d_i)$, $\SO_{\Bbb P^1}(r_i\alpha_i)$ are divisors of degree $d_i$,
$r_i\alpha_i$ on $C$, $\Bbb P^1$ respectively.
Here we remark that for any torsion free sheaf $F_i$ on $X$ we have
$c_1(F_i)|_{f^{-1}(\bullet)}=c_1(F_i|_{f^{-1}(\bullet)})$ for general 
points on $C$ 
(resp. $c_1(F_i)|_{\pi^{-1}(\bullet)}=c_1(F_i|_{\pi^{-1}(\bullet)})$ 
for general points on $\Bbb P^1$). Therefore
$d_i$ are the degrees of $F_i$ on the general fiber of $\pi$ respectively.
Without confusion, we denote the degree of $F_i$ (resp. $E_i$) on the 
generic fiber of $\pi$ by $deg(F_i)$ (resp. $deg(E_i)$).
Consequently, $\mu(E_i)$, $\mu(E)$ denote the slope of restrictions
of $E_i$, $E$ to the generic fiber of $\pi$ respectively. Note that
$$c_2(F_i')=c_2(F_i)-(r_i-1)c_1(F_i)\pi^*\SO_{\Bbb P^1}(\alpha_i)=
c_2(F_i)-(r_i-1)d_i\alpha_i,$$
$c_1(F_i)^2=2r_id_i\alpha_i$ and $c_1(E)^2=2d\sum^n_{i=1}r_i\alpha_i$,
we have
$$\Delta(E)=2r\left(\sum_{i=1}^nc_2(F'_i)+\mu(E)\sum_{i=1}^nr_i\alpha_i
-\sum^n_{i=1}d_i\alpha_i\right).$$
Let $rk(E_i)$ denote the rank of $E_i$, note that $r_i=rk(E_i)-rk(E_{i-1})$
and $d_i=deg(E_i)-deg(E_{i-1})$, we have
$$\Delta(E)=2r\left(\sum_{i=1}^nc_2(F'_i)+\sum_{i=1}^{n-1}
(\mu(E)-\mu(E_i))(\alpha_i-\alpha_{i+1})rk(E_i)
\right).\tag2.2$$

\proclaim{Lemma 2.1} Any torsion free sheaf $\SE$ of rank $r$ 
on a ruled surface, with
generic splitting type $(0^{\oplus r})$, must have $c_2(\SE)\ge 0$.
\endproclaim

\demo{Proof} The argument is in fact contained in the proof of
Lemma 1.4 of [GL]. One can choose a divisor $D$ on $C$ such that
$f^*\SO_C(D)\subset \SE$ and $\SE/f^*\SO_C(D)$ is torsion free.
Since $\SE/f^*\SO_C(D)$ has generic splitting type $(0^{\oplus (r-1)})$,
by induction hypothesis on rank, we can assume that
$c_2(\SE/f^*\SO_C(D))\ge 0$. Hence
$$c_2(\SE)=c_2(\SE/f^*\SO_C(D))+f^{-1}(D)(c_1(\SE)-f^{-1}(D))
 =c_2(\SE/f^*\SO_C(D))\ge 0.$$
\enddemo

\proclaim{Proposition 2.2} If $g\ge 4$, then $\Delta(E)\ge 2r$.
The equality holds if and only if $E$ has generic splitting type
$(0^{\oplus r})$ and $c_2(E)=1$.\endproclaim

\demo{Proof} Since the rational curve passes through a generic point 
$[W]\in M$, then the bundle $E|_{\pi^{-1}([W])}\cong W$ is $(1,1)$-stable.
Thus, using the inequality (2.2), we have
$$ \Delta(E)\ge 2r\left(\sum_{i=1}^nc_2(F'_i)+\sum_{i=1}^{n-1}
(\alpha_i-\alpha_{i+1})\right),$$
the equality holds if and only if $n=1$. 

\newpage

By Lemma 2.1, if $n>1$, we have
$$\Delta(E)>2r\sum_{i=1}^{n-1}
(\alpha_i-\alpha_{i+1})\ge 2r.$$
If $n=1$, then $\Delta(E)=2rc_2(E)\ge 2r$ 
(note that we have assumed $\alpha_n=0$). The equality $\Delta(E)=2r$ holds
if and only if $c_2(E)=1$.
\enddemo

From now on, we assume that $E$ has generic splitting type $(0^{\oplus r})$.
If $E$ has a jumping line $X_p=f^{-1}(p)$ ($p\in C$), i.e.,
$$E|_{X_p}=\bigoplus^n_{i=1}\SO_{X_p}(\beta_i)^{\oplus r_i},
\quad \beta_1>\,\cdots\,>\beta_n$$
with the type $(\beta_1^{\oplus r_1},\,\cdots\,,\beta_n^{\oplus r_n})$
different from $(0^{\oplus r})$. Then we can perform the elementary
transformation on $E$ along $X_p$ by taking $F$ to be the kernel of the
(unique surjective) homomorphism
$\phi: E\to E|_{X_p}\to \SO_{X_p}(\beta_n)^{\oplus r_n}$. Clearly,
$$0@>>>F@>>>E@>\phi>>\SO_{X_p}(\beta_n)^{\oplus r_n}@>>>0.\tag 2.3$$
An easy calculation yields

\proclaim{Lemma 2.3} $c_1(F)=c_1(E)-r_nX_p$ and $c_2(F)=c_2(E)+r_n\beta_n$.
\endproclaim

\demo{Proof} By the exact sequence (2.3), the computation is 
straightforward.
\enddemo

\proclaim{Lemma 2.4} 
If $c_2(E)=1$ and $E$ has generic splitting type $(0^{\oplus r})$, then
$E$ has exactly one jumping line $X_p$ and the elementary
transformation $F$ along $X_p$ is isomorphic to $f^*V$ for a vector bundle
$V$ over $C$.
\endproclaim

\demo{Proof} The $E$ has at least one jumping line. Otherwise, $E$
will be a pullback of a vector bundle over $C$, which is impossible.
At any jumping line $X_p$, with splitting type
$(\beta_1^{\oplus r_1},\,\cdots\,,\beta_n^{\oplus r_n})$, we must have
$\beta_n<0$. Hence, by Lemma 2.3 and Lemma 2.1, $E$ has a unique
jumping line $X_p$ with $\beta_n=-1$ and $r_n=1$. Then $F$ has no jumping
line, thus $F=f^*V$ for a vector bundle $V$ over $C$.   
\enddemo

Therefore by Proposition 2.2 and Lemma 2.4, if $\Delta(E)=2r$, we have
$$0@>>>f^*V@>>>E@>\phi>>\SO_{X_p}(-1)@>>>0.\tag2.4$$

\proclaim{Proposition 2.5} If $g\ge 4$ and $\Delta(E)=2r$, then the 
rational curve
is a Hecke curve.\endproclaim

\demo{Proof}  For any $x\in\Bbb P^1$, 
let $E_x$ denote $E|_{C\times\{x\}}$. 
Restrict the sequence (2.4) to
$\pi^{-1}(x)=C\times\{x\}$, we get
$$ 0@>>>V@>>>E_x@>\phi_x>>\SO_{X_p}(-1)_x@>>>0.\tag2.5$$
Since the curve passes through the generic point $[W]\in M$, there is a
$x_0\in \Bbb P^1$
such that $E_{x_0}\cong W$ is
$(1,1)$-stable by Lemma 1.1. Then $V$ is a $(1,0)$-stable bundle by
Lemma 1.2 (thus when $x=x_0$, the sequence (2.5) corresponds to the
sequence (1.1) and $V$ corresponds to $W^{\zeta}$ in the definition of
Hecke curves, see {\S1}). 

\newpage
Let
$\iota_x: V_p\to E_x|_p=E_{(p,x)}$
be the homomorphism between the fibers at $p$ induced by the sheaf injection
$V\to E_x$ 
in sequence (2.5). Then the kernel $ker(\iota_x)$ is a $1$-dimensional
subspace of $V_p$. 
When $x$ moves on $\Bbb P^1$, these 
$[ker(\iota_x)]\in\Bbb P(V_p)$ form a line $\SH\subset\Bbb P(V_p)$.
It is easy to check that, as the same as (1.3),  for any $x\in\Bbb P^1$
$${\wt W}^{ker(\iota_x)}\cong E_x^{\vee}.$$
Thus $\{({\wt W}^{ker(\iota_x)})^{\vee};\,\,[ker(\iota_x)]\in\SH\}$ defines
the given rational curve. That is, the given rational curve is a Hecke curve.
\enddemo

\proclaim{Theorem 2.6} Assume that $g\ge 4$. Then any rational curve
of $M$ passing through the generic point of $M$ has at least
degree $2r$ with respect to $-K_M$. It has degree $2r$ if and only if
it is a Hecke curve.\endproclaim

\demo{Proof} By (2.1), the degree $-K_{M}\cdot \Bbb P^1$ equals to
$\Delta(E)$. Then, by Proposition 2.2, it has degree at least $2r$.
If it has degree $2r$, then by Proposition 2.5 it must be a Hecke curve.
It was known that any Hecke curve has degree $2r$. We are done
\enddemo

\enddocument